\documentclass[a4paper
,12pt
,leqno
,twoside
]{article}

\usepackage{amsmath, amsthm, amssymb,bm}
\numberwithin{equation}{section}

\newcommand{\bC}{\mathbb{C}}

\newcommand{\bE}{\mathbb{E}}

\newcommand{\mr}[1]{\mathrm{#1}}
\newcommand{\mcal}[1]{\mathcal{#1}}

\def\({ \left( }
\def\){ \right)}

\def\trans#1{\mathord{\mathopen{{\vphantom{#1}}^t}#1}}

\theoremstyle{plain}
\newtheorem{thm}{Theorem}
\newtheorem{prop}{Proposition}
\newtheorem{lem}[prop]{Lemma}
\newtheorem{cor}[thm]{Corollary}
\theoremstyle{definition}

\newtheorem{example}{Example}
\newtheorem{remark}{Remark}
\theoremstyle{conjecture}

\theoremstyle{problem}

\title{\bfseries General moments of matrix
elements from circular orthogonal ensembles}
\author{\textsc{Sho Matsumoto} 
\\
Graduate School of Mathematics, Nagoya University, \\
 Nagoya 464-8602, Japan. \\
\qquad E-mail: sho-matsumoto@math.nagoya-u.ac.jp
}

\date{\empty}

\begin{document}

\maketitle

\begin{abstract}
The aim of this paper is to present
a systematic method for computing moments of matrix elements taken 
from circular orthogonal ensembles (COE).
The formula is given as a sum of Weingarten functions for orthogonal groups
but the technique for its proof involves Weingarten calculus for unitary groups.
As an application, explicit expressions for the moments of a single matrix element
of a COE matrix are given.

\end{abstract}

\section{Introduction and main result}

The study of distributions for matrix elements of
random matrices 
is one of main themes in random matrix theory.
When a random matrix $X=(x_{ij})_{1 \le i,j \le N}$ is given,
we would like to know how to compute the moments
$$
\bE[x_{i_1 j_1}  x_{i_2 j_2} \cdots x_{i_n j_n}] \qquad 
\text{or} \qquad 
\bE[x_{i_1 j_1}  x_{i_2 j_2} \cdots x_{i_n j_n}
\overline{x_{k_1 l_1}  x_{k_2 l_2} \cdots x_{k_m l_m}}].
$$
There is a large body of research in this area, but we shall focus on recent papers:
\begin{itemize}
\item Gaussian matrix: There are three typical Gaussian matrix ensembles:
Gaussian orthogonal, unitary, and symplectic ensembles (GOE, GUE, and GSE);
see, e.g., \cite{Mehta}.
The calculus of moments for these is well known, under the name Wick calculus;
see, e.g., \cite{Zv}.
\item Haar-distributed unitary matrix: 
The compact Lie group $U(N)$ of $N \times N$ unitary matrices has 
the Haar probability measure. 
Collins \cite{C} showed how to compute the moments and called his technique
Weingarten calculus, a reference to the trailblazing work of Weingarten \cite{W}.
The main tool in this theory is the unitary Weingarten function $\mathrm{Wg}^{U(N)}_n$,
which is a class function on the symmetric group $S_{n}$. See also \cite{Samuel, MN}.
\item Haar-distributed orthogonal matrix: The compact Lie group
$O(N)$ of $N \times N$ real orthogonal matrices also has the Haar probability measure.
Weingarten calculus for $O(N)$ has been developed in \cite{CS, CM, Mat_JMortho}.
In this calculus, the necessary function is the orthogonal Weingarten function $\mathrm{Wg}^{O(N)}_n$,
which is a function on $S_{2n}$ with an invariant property under the two-sided action of
the hyperoctahedral group.
\item Central Wishart matrix : A Wishart random matrix is a matrix version of 
a chi-square distribution.
Graczyk, Letac, and Massam \cite{GLM1, GLM2} established the moment computations for 
the complex and real cases.
\item Inverse of a central Wishart matrix: The distribution of the inverse $W^{-1}$ of $W$
is quite different from the original Wishart matrix $W$.
The inverse distributions of the complex and real cases were studied in \cite{GLM1} and \cite{Mat_Wishart},
respectively. 
For reasons not fully understood, unitary and orthogonal Weingarten functions are necessary
for the complex and real cases, respectively.
\item Noncentral Wishart matrix: The noncentral version of the Wishart matrix
is derived from noncentral Gaussian vectors. Kuriki and Numata \cite{KN} showed 
computation of the moments using matchings.
Note that the problem for the inverted  noncentral Wishart matrix is still open.
\end{itemize}

There are three much-studied circular ensembles of random unitary matrices,
introduced by Dyson:
circular orthogonal, unitary, and symplectic ensembles (COE, CUE, and CSE).
The density function for their eigenvalues $\Lambda_1,\Lambda_2,\dots,\Lambda_N$
is proportional to $\prod_{1 \le i<j \le N} |\Lambda_i-\Lambda_j|^\beta$
with $\beta=1,2,4$ for COE, CUE and CSE, respectively.
See \cite{Mehta} for details.

A CUE matrix, that is, a random matrix taken from
the CUE, is simply a Haar-distributed unitary matrix from $U(N)$,
and such a matrix has been studied as stated above. 
Therefore, the next obvious object of study is the COE.
The COE is the probability space of symmetric unitary matrices
defined by the property of being invariant under 
automorphisms
$V \to \trans{U_0} V U_0$,
where $U_0$ is a fixed unitary matrix.
The COE is a realization of the compact symmetric space $U(N)/O(N)$
equipped with the probability measure derived from the 
Haar measure on $U(N)$.

Our aim in the present paper is to establish the calculus for
the moments of matrix elements of COE matrices.
If $U$ is an $N \times N$ Haar-distributed unitary matrix,
then an $N \times N$ COE matrix $V$ can be written as
$V= \trans{U} U$. 
In other words,
if $V=(v_{ij})$ and $U=(u_{ij})$
then
$v_{ij}= \sum_{k=1}^N u_{ki} u_{kj}$.
Therefore, the desired moment for $V$
can be written as a sum of moments for $U$.
Since the computation of moments for $U$ is known,
the moments for $V$ can in principle be computed.
(In this direction, Collins and Stolz \cite{CSt} studied 
asymptotic behaviors for the COE and for random matrices 
associated with classical compact symmetric spaces.)
Herein, we will express the COE moment 
in terms of unitary Weingarten functions $\mathrm{Wg}^{U(N)}_{2n}$ (Proposition \ref{prop:eq:M_unitary}).
However, this expression is ill-suited to practical computation,
so obtaining this expression is not our main goal in the present paper. 

Our main result is given as follows.
Put  $[N]=\{1,2,\dots,N\}$ and define the right action of the symmetric group $S_{2n}$ on $[N]^{2n}$ by
$$
\bm{i}^\sigma =(i_{\sigma(1)},i_{\sigma(2)},\dots, i_{\sigma(2n)}) \qquad
\text{for $\bm{i}=(i_1,i_2,\dots,i_{2n}) \in [N]^{2n}$ and $\sigma \in S_{2n}$}.
$$

\begin{thm} \label{thm:main}
Let $N$ be a positive integer, and 
let $V=(v_{ij})_{1 \le i,j \le N}$ be an $N \times N$ COE matrix.
Let $\bm{i}=(i_1,\dots,i_{2n})$ and 
$\bm{j}=(j_1,\dots,j_{2n})$ be two
sequences in $[N]^{2n}$, and 
let
$$
M_N(\bm{i},\bm{j})= \bE[v_{i_1, i_2} v_{i_3, i_4} \cdots v_{i_{2n-1},i_{2n}}
\overline{v_{j_1, j_2} v_{j_3,j_4} \cdots v_{j_{2n-1},j_{2n}}}]
$$
and 
$$
\mcal{M}(\bm{i},\bm{j};z)= \sum_{\begin{subarray}{c} \sigma \in S_{2n} \\
\bm{j}=\bm{i}^\sigma \end{subarray}}
\mathrm{Wg}^{O}_n (\sigma;z) \qquad (z \in \bC),
$$
where the sum runs over all permutations $\sigma \in S_{2n}$ satisfying $\bm{j}=\bm{i}^\sigma$,
and $\mathrm{Wg}^{O}_n(\sigma;z)$ is the orthogonal Weingarten function  
defined 
in \eqref{eq:DefWg} below.
Then we have
$$
M_N(\bm{i},\bm{j})= \mcal{M}(\bm{i},\bm{j};N+1).
$$
\end{thm}

If $N +1 \ge n$, our formula can be written as 
$$
M_N(\bm{i},\bm{j})= 
\sum_{\begin{subarray}{c} \sigma \in S_{2n} \\
\bm{j}=\bm{i}^\sigma \end{subarray}}
\mathrm{Wg}^{O(N+1)}_n (\sigma)
$$
with $\mathrm{Wg}^{O(N+1)}_n(\sigma):= \mathrm{Wg}^O_n(\sigma;N+1)$.
However,
if $N+1 <n$ and $\sigma \in S_{2n}$, then
the complex function $z \mapsto \mathrm{Wg}_n^{O}(\sigma;z)$ may have a pole at 
$z=N+1$, and so $\mathrm{Wg}^O_n(\sigma;N+1)$ does not make sense.
Nevertheless, $\mcal{M}(\bm{i},\bm{j};N+1)$ is well defined after some cancellations of poles.

\begin{example}
Let us compute $\bE[|v_{11}|^6]=
\bE[v_{11}v_{11}v_{11} \overline{v_{11}v_{11}v_{11}}]=M_N(\bm{i},\bm{i})$
with $\bm{i}=(1,1,1,1,1,1)$.
Theorem \ref{thm:main} implies 
$M_N(\bm{i},\bm{i}) = [\sum_{\sigma \in S_{6}} \mathrm{Wg}_3^{O}(\sigma;z)]_{z=N+1}$.
A representation-theoretic discussion (the proof of Theorem \ref{thm:diagonal}) gives
$\sum_{\sigma \in S_6} \mathrm{Wg}^O_3(\sigma;z) = \frac{48}{z(z+2)(z+4)}$;
hence,
$$
\bE[|v_{11}|^6]= \frac{48}{(N+1)(N+3)(N+5)}.
$$
If $N+1 <n=3$, i.e., if $N=1$, we cannot replace
$[\sum_{\sigma \in S_{6}} \mathrm{Wg}_3^{O}(\sigma;z)]_{z=N+1=2}$
by $\sum_{\sigma \in S_{6}} \mathrm{Wg}_3^{O}(\sigma;2)$
since
$\mathrm{Wg}^O_3(\mathrm{id}_6;z)=\frac{z^2+3z-2}{z(z+2)(z+4)(z-1)(z-2)}$
has a pole at $z=2$.
\end{example}

Surprisingly, the orthogonal Weingarten function has appeared three times now
in random matrix theory:
the first appearance was, just as its name implies, in the study of Haar-distributed orthogonal matrices \cite{CM} 
(implicitly \cite{CS})
and the second was related to inverted central real Wishart matrices \cite{Mat_Wishart}.
Furthermore,
the reader should note that 
the function $\mathrm{Wg}^{O(N+1)}_n$ appears, not $\mathrm{Wg}^{O(N)}_n$,
in our formula.
The formula is derived by using Weingarten calculus for $U(N)$ and 
harmonic analysis for symmetric groups, and the function $\mathrm{Wg}^{O(N+1)}_n$ 
then appears 
as a result. 
Since the COE is identified with the compact symmetric space $U(N)/O(N)$, 
it is small wonder that 
an orthogonal Weingarten function appears.
However, 
the theoretical reason why we need $O(N+1)$ rather than $O(N)$
seems quite mysterious. 

The present paper is organized as follows.
In section \ref{sec:Wg}, we review the definition of orthogonal Weingarten functions 
developed in \cite{CM, Mat_JMortho, Mat_Wishart}.
In section \ref{section:Proof}, we prove Theorem \ref{thm:main}.
In section \ref{section:single moment}, 
using Theorem \ref{thm:main},
we derive explicit expressions 
for the moments $\bE[|v_{ij}|^{2n}]$ of a single matrix element $v_{ij}$.
In section \ref{section:asymptotic},
we show that there exist nonnegative integers $a$ and $b$ satisfying
$$
M_N(\bm{i},\bm{j})= a N^{-n} -b N^{-n-1} +O(N^{-n-2}), \qquad N \to \infty,
$$
and evaluate $a,b$.
In the final short section, section \ref{section:example}, we show 
how to use Theorem \ref{thm:main} for cases where 
the degree of integrand $n$ is small.

Theorem \ref{thm:main} is a very powerful tool for computing various averages on the COE in the same way that Wick calculus for Gaussian ensembles and Weingarten calculus for classical groups
are.
We can expect many applications of Theorem \ref{thm:main} to appear 
in future research.

\begin{remark}
In a recent preprint \cite{M_COE}, 
the author showed obtaining some results for simple cases of $M_N(\bm{i},\bm{j})$
by the same method as in the present paper.
Results in \cite{M_COE} seem unsatisfactory
compared to the present ones.
All these previous results are
included in the present note with much simpler proofs.
\end{remark}

\begin{remark}
The COE is a random matrix ensemble associated with 
the compact symmetric space $U(N)/O(N)$,
which is of Cartan class AI.
Other random matrices from classical compact symmetric spaces
are expressed in terms of Haar-distributed matrices from classical Lie groups.
Thus, if we apply Weingarten calculus for classical groups to these,
we obtain the same formulas as those given in Theorem \ref{thm:main}.
We leave them as future research.
\end{remark}

\medskip

\noindent
{\it Acknowledgment.}
The author would like to acknowledge an inspiring conversations with Benoit Collins.
This work was supported by a Grant-in-Aid for Young Scientists (B) No. 22740060.

\section{Orthogonal Weingarten functions}  \label{sec:Wg}

We first review the theory of Weingarten functions for orthogonal groups;
see \cite{CM,Mat_JMortho, Mat_Wishart} for details.
Claims in the first two subsections can be found in \cite[VII.2]{Mac}.

\subsection{Hyperoctahedral groups and coset-types} \label{subsec:Hyper}

Let $H_n$ be the subgroup in $S_{2n}$ generated by
transpositions $(2k-1 \ 2k)$, $1 \le k \le n$, and 
by double transpositions $(2i-1 \  2j-1) (2i \ 2j),$ \ $1 \le i<j \le n$.
The group $H_n$ is called the hyperoctahedral group.
Note that $|H_n|=2^n n!$, and the pair $(S_{2n}, H_n)$ is
a Gelfand pair in the sense of \cite[VII]{Mac}.

A partition $\lambda=(\lambda_1,\lambda_2,\dots)$ is a weakly decreasing sequence
of nonnegative integers such that $|\lambda|:=\sum_{i \ge 1} \lambda_i$ is finite.
If $|\lambda|=n$, we call $\lambda$ a partition of $n$ and write $\lambda \vdash n$.
The length $\ell(\lambda)$ of $\lambda$ is defined as the number of nonzero $\lambda_i$.

Given $\sigma \in S_{2n}$, we attach an undirected graph $\Gamma(\sigma)$ with vertices
$1,2,\dots,2n $ and edge set consisting of 
$\big\{ \{2k-1,2k\} \ | \ k \in [n] \big\}$ and 
$\big\{ \{ \sigma(2k-1), \sigma(2k)\} \ | \ k \in [n] \big\}$.
Each vertex of the graph lies on exactly two edges, and
the number of vertices in each connected component is even.
If the numbers of vertices are
$2\lambda_1 \ge 2\lambda_2 \ge \dots \ge 2\lambda_l$, 
then the sequence $\lambda=(\lambda_1,\lambda_2,\dots,\lambda_l)$ is a partition of $n$.
We will refer to $\lambda$ as the coset-type of $\sigma$.

In general, given $\sigma,\tau \in S_{2n}$, their coset-types coincide if and only if
$H_n \sigma H_n= H_n \tau H_n$.
Hence, we have the double coset decomposition of $H_n$ in $S_{2n}$
$$
S_{2n}= \bigsqcup_{\mu \vdash n} H_\mu, \qquad \text{where $
H_\mu= \{\sigma \in S_{2n} \ | \ \text{the coset-type of $\sigma$ is $\mu$}\} $}.
$$
Note that $H_{(1^n)}=H_n$ and 
$|H_\mu|= (2^n n!)^2/ (2^{\ell(\mu)} z_\mu)$, where
$$
z_\mu= \prod_{r \ge 1} r^{m_r(\mu)} m_r(\mu)!
$$
with multiplicities $m_r(\mu)= | \{ i \ge 1 \ | \ \mu_i=r\}|$.

Denote by $\ell'(\sigma)$ the number of connected components of $\Gamma(\sigma)$,
or equivalently, $\ell'(\sigma)=\ell(\mu)$ if $\sigma \in H_\mu$.

\subsection{Zonal spherical functions} \label{subsec:zonal}

For two complex-valued functions $f_1,f_2$ on $S_{2n}$,
their convolution $f_1 * f_2$ is defined by
$$
(f_1* f_2)(\sigma)= \sum_{\tau\in S_{2n}} f_1(\sigma \tau^{-1}) f_2(\tau) \qquad (\sigma \in S_{2n}).
$$
Let $e_{H_n}$ be the function on $S_{2n}$ defined by
$$
e_{H_n}(\sigma) = 
\begin{cases}
(2^n n!)^{-1} & \text{if $\sigma \in H_n$,} \\
0 & \text{otherwise}.
\end{cases}
$$
For each $\lambda \vdash n$, we define
the zonal spherical function $\omega^\lambda$ by
$\omega^\lambda=\chi^{2\lambda} * e_{H_n}= e_{H_n} *\chi^{2\lambda}$, or by
$$
\omega^\lambda(\sigma)= \frac{1}{2^n n!} \sum_{\zeta \in H_n} \chi^{2\lambda}(\sigma \zeta)
\qquad
(\sigma \in S_{2n}),
$$
where $\chi^{2\lambda}$ is the irreducible character of $S_{2n}$ associated with
$2\lambda=(2\lambda_1,2\lambda_2,\dots)$.
The functions $\omega^\lambda$ are $H_n$-biinvariant,
and so they take a constant value on each $H_\mu$.

We will use the following two lemmas later.

\begin{lem} \label{lem:omega-chi}
If $\lambda \vdash n$ and $\mu \vdash 2n$, then
$\omega^\lambda * \chi^\mu= \frac{(2n)!}{f^{2\lambda}} \delta_{2\lambda,\mu} \omega^\lambda$.
Here, $f^{2\lambda}$ is the dimension of the irreducible representation of character $\chi^{2\lambda}$. Hence, 
$\omega^\lambda * \omega^\rho= \frac{(2n)!}{f^{2\lambda}} \delta_{\lambda,\rho} \omega^\lambda$
for $\lambda,\rho \vdash n$.
\end{lem}
\begin{proof}
The well-known fact 
$\chi^\mu* \chi^\rho= \frac{(2n)!}{f^\mu} \delta_{\mu,\rho} \chi^\mu$ gives
$\omega^\lambda*\chi^\mu= 
e_{H_n}* \chi^{2\lambda} *\chi^\mu= 
\frac{(2n)!}{f^{2\lambda}} \delta_{\mu,2\lambda} e_{H_n} * \chi^{2\lambda}
=\frac{(2n)!}{f^{2\lambda}} \delta_{\mu,2\lambda} \omega^\lambda$.
\end{proof}

\begin{lem}
For any complex number $z$ and $\sigma \in S_{2n}$,
\begin{equation} \label{eq:ExpandN}
z^{\ell'(\sigma)}= \frac{2^n n!}{(2n)!} \sum_{\lambda \vdash n} f^{2\lambda} C'_\lambda(z) 
\omega^\lambda (\sigma) 
\end{equation}
with
\begin{equation} \label{eq:defC'}
C'_\lambda(z) =\prod_{(i,j) \in \lambda} (z+2j-i-1) = \prod_{i=1}^{\ell(\lambda)}
\prod_{j=1}^{\lambda_i} (z+2j-i-1).
\end{equation}
\end{lem}

\begin{proof}
This is a special case of an identity involving power-sum symmetric polynomials
and zonal polynomials.
See \cite[(4.9)]{Mat_Wishart}.
\end{proof}

\subsection{Orthogonal Weingarten functions} \label{subsec:Wg}

For each positive integer $n$,
we define the {\it orthogonal Weingarten function} as
\begin{equation} \label{eq:DefWg}
\mr{Wg}^{O}_n(\sigma;z) = \frac{2^n n!}{(2n)!} 
\sum_{\lambda \vdash n}
\frac{f^{2\lambda}}{C_\lambda'(z)} 
\omega^\lambda(\sigma) \qquad (\sigma \in S_{2n}, \ z \in \bC),
\end{equation}
which was first defined in \cite{CM}.
If we fix $\sigma$, the complex function $z \mapsto \mathrm{Wg}^O_n(\sigma;z)$
is a rational function in $z$. 
Since $C'_\lambda(z)$ has a factor $(z-1)(z-2)\cdots (z-\ell(\lambda)+1)$,
the function $\mr{Wg}^{O}_n(\sigma;z)$ may have poles at the points $z=1,2,\dots,n-1$ 
on the half line $\mathbb{R}_{>0}$.

If $N  \ge n$, then $C_\lambda'(N) \not=0$ for all $\lambda \vdash n$, and 
we can  define
$$
\mathrm{Wg}^{O(N)}_n(\sigma)= [\mathrm{Wg}^{O}_n(\sigma;z)]_{z=N}=
\frac{2^n n!}{(2n)!} 
\sum_{\lambda \vdash n}
\frac{f^{2\lambda}}{C_\lambda'(N)} 
\omega^\lambda(\sigma) \qquad (\sigma \in S_{2n}).
$$
We will refer to the above as the {\it Weingarten function for the orthogonal group $O(N)$}.

The function $S_{2n} \ni \sigma \mapsto \mathrm{Wg}^O_n(\sigma;z)$ is constant on each
double coset $H_\mu$.
We will denote by $\mathrm{Wg}^O_n(\mu;z)$ the value 
of $\mathrm{Wg}^O_n(\sigma;z)$
at $\sigma \in H_\mu$.

\begin{example} \label{ex:WgO}
The following examples can be found in \cite{CM,CS}.
$$
\mathrm{Wg}^O_1([1];z)= \frac{1}{z},
$$
$$
\mathrm{Wg}^O_2([2];z)= \frac{-1}{z(z+2)(z-1)}, \qquad 
\mathrm{Wg}^O_2([1^2];z)= \frac{z+1}{z(z+2)(z-1)},
$$
\begin{align*}
\mathrm{Wg}^O_3([3];z)=& \frac{2}{z(z+2)(z+4)(z-1)(z-2)}, \\
\mathrm{Wg}^O_3([2,1];z)=& \frac{-1}{z(z+4)(z-1)(z-2)}, \\
\mathrm{Wg}^O_3([1^3];z)=& \frac{z^2+3z-2}{z(z+2)(z+4)(z-1)(z-2)}.
\end{align*}
\end{example}

\begin{remark}
The orthogonal Weingarten function was introduced in \cite{CS} implicitly,
in a study of moments for Haar-distributed orthogonal matrices. 
The explicit definition  was first given in \cite{CM}.
If $O=(o_{ij})_{1 \le i,j \le N}$ is a Haar-distributed orthogonal matrix taken from $O(N)$
and $N \ge n$,
then for each $\sigma \in S_{2n}$,
$$
\mathrm{Wg}^{O(N)}_n(\sigma)= \bE [o_{1, j_1}
o_{1, j_2} o_{2, j_3} o_{2, j_4} \cdots o_{n, j_{2n-1}} o_{n, j_{2n}}] 
$$
with $(j_1,j_2,\dots,j_{2n}) =(1,1,2,2\dots,n,n)^{\sigma}$.
\end{remark}

\section{Proof of Theorem \ref{thm:main}} \label{section:Proof}

\subsection{Haar-distributed unitary matrix}

The orthogonal Weingarten function appears in the statement of Theorem \ref{thm:main},
but we must employ a unitary Weingarten function in our proof.

We define the {\it unitary Weingarten function} as
\begin{equation} \label{eq:unitaryWg}
\mr{Wg}^{U}_n(\sigma;z) = \frac{1}{n!}
 \sum_{\lambda \vdash n } 
\frac{f^\lambda}{C_\lambda(z)} \chi^\lambda(\sigma)
\qquad (\sigma \in S_n, \ z \in \mathbb{C}).
\end{equation}
Here, $\chi^\lambda$ is the irreducible character of $S_n$ associated with $\lambda$,
$f^\lambda$ is its degree, 
and 
$$
C_\lambda(z)=\prod_{(i,j)\in\lambda} (z+j-i).
$$
As in the case of the orthogonal Weingarten function,
the function $\mr{Wg}^{U}_n(\sigma;z)$ may have poles at the points $z=1,2,\dots,n-1$ 
on $\mathbb{R}_{>0}$.

If $N  \ge n$, we let
$$
\mathrm{Wg}^{U(N)}_n(\sigma)= [\mathrm{Wg}^{U}_n(\sigma;z)]_{z=N}=
\frac{1}{n!}
 \sum_{\lambda \vdash n } 
\frac{f^\lambda}{C_\lambda(N)} \chi^\lambda(\sigma) \qquad (\sigma \in S_{n})
$$
and refer to it as the {\it Weingarten function for the unitary group $U(N)$}.

Given a unitary matrix $U =(u_{ij})\in U(N)$ and 
four sequences 
$\bm{i}=(i_1,i_2,\dots,i_n)$, $\bm{j}=(j_1,\dots,j_n)$, 
$\bm{i}'=(i_1',\dots,i_m')$,
$\bm{j}'=(j_1,\dots,j_m')$ of indices in $[N]$,
we let
$$
U(\bm{i},\bm{j}|\bm{i}',\bm{j}')=
u_{i_1 j_1} u_{i_2 j_2} \cdots u_{i_n j_n} 
\overline{u_{i_1' j_1'} u_{i_2' j_2'} \cdots u_{i_m' j_m'}}.
$$

\begin{lem}[Weingarten formula for the unitary group; see, e.g., \cite{C,CS}] \label{lem:WgFormula}
Let $U$ be an $N \times N$ Haar-distributed unitary matrix and let
$\bm{i}=(i_1,i_2,\dots,i_n)$, $\bm{j}=(j_1,\dots,j_n)$, 
$\bm{i}'=(i_1',\dots,i_m')$,
$\bm{j}'=(j_1,\dots,j_m')$ be four sequences of indices in $[N]$.
If $n=m$, 
$$
\bE[U(\bm{i},\bm{j}|\bm{i}',\bm{j}')]= \Bigg[
\sum_{\begin{subarray}{c} \sigma \in S_n \\ \bm{i}^\sigma=\bm{i}' \end{subarray}} 
\sum_{\begin{subarray}{c} \tau \in S_n \\ \bm{j}^\tau=\bm{j}' \end{subarray}} \mr{Wg}^{U}_n(\sigma^{-1} \tau;z) \Bigg]_{z=N},
$$
and it vanishes otherwise.
\end{lem}

\begin{example}
If $U=(u_{ij})$ is an $N \times N$ Haar-distributed unitary matrix, then
\begin{align*}
&\bE[|u_{11}|^4] = \Big[\sum_{\sigma \in S_2} \sum_{\tau \in S_2} \mathrm{Wg}^{U}_2
(\sigma^{-1} \tau;z) \Big]_{z=N} 
= [2 \mathrm{Wg}^U_2(\mathrm{id}_2;z)+ 2 \mathrm{Wg}^U_2( (1 \ 2);z)]_{z=N} \\ 
=& \Big[2 \frac{1}{(z+1)(z-1)} + 2 \frac{-1}{z(z+1)(z-1)}\Big]_{z=N} 
= \Big[\frac{2}{z(z+1)} \Big]_{z=N}= \frac{2}{N(N+1)}
\end{align*}
for all $N \ge 1$.
\end{example}

\subsection{Expression in terms of unitary Weingarten functions}

Let $V=(v_{ij})_{1 \le i,j \le N}$ be an $N \times N$ COE matrix.
Let $M(\pi)=(\delta_{i,\pi(j)})_{1 \le i,j \le N}$ be the permutation matrix
associated to $\pi \in S_{N}$.
Since $\trans{M(\pi)} V M(\pi)$ has the same distribution as $V$,
we have
$$
v_{ii} \stackrel{\mathrm{d}}{=} v_{11}, \qquad 
v_{ij} \stackrel{\mathrm{d}}{=} v_{12} \quad (i \not=j),
$$
where $X\stackrel{\mathrm{d}}{=} Y$ means that
$X$ has the same distribution as $Y$.

Using a Haar-distributed unitary matrix $U=(u_{ij})_{1 \le i,j \le N}$,
we can write $V=\trans{U} U$.
Let  $\bm{i}=(i_1,i_2,\dots,i_{2n}) \in [N]^{2n}$ and 
$\bm{j}=(j_1,j_2,\dots,j_{2m}) \in [N]^{2m}$. 
Since $v_{ij}=\sum_k u_{ki} u_{kj}$,
we have
\begin{equation} \label{eq:v_u}
\bE[v_{i_1, i_2} v_{i_3, i_4} \cdots v_{i_{2n-1}, i_{2n}} 
\overline{v_{j_1, j_2} v_{j_3, j_4} \cdots v_{j_{2m-1}, j_{2m}} }]
=\sum_{\bm{k} \in [N]^n} \sum_{\bm{l} \in [N]^m}
\bE[U(\tilde{\bm{k}}, \bm{i}| \tilde{\bm{l}},\bm{j})],
\end{equation}
where 
\begin{align*}
\bm{k}=& (k_1,k_2,\dots,k_n), & \bm{l}=& (l_1,l_2,\dots,l_m), \\
\tilde{\bm{k}}=&(k_1,k_1,k_2,k_2,\dots,k_n,k_n), &  
\tilde{\bm{l}}=&(l_1,l_1,l_2,l_2,\dots,l_m,l_m).
\end{align*}
Using \eqref{eq:v_u} and Lemma \ref{lem:WgFormula} gives the following.

\begin{lem} \label{lem:COEcondition}
$\bE[v_{i_1, i_2} v_{i_3, i_4} \cdots v_{i_{2n-1}, i_{2n}} 
\overline{v_{j_1, j_2} v_{j_3, j_4} \cdots v_{j_{2m-1}, j_{2m}} }]$ vanishes
unless $n=m$ and there exists a permutation $\sigma \in S_{2n}$
such that $\bm{j}=\bm{i}^\sigma$.
\end{lem}

Hereafter, we will assume that 
$\bm{i}=(i_1,\dots,i_{2n})$ and 
$\bm{j}=(j_1,\dots,j_{2n})$ and let 
$$
M_N(\bm{i},\bm{j})
=\bE[v_{i_1, i_2} v_{i_3, i_4} \cdots v_{i_{2n-1}, i_{2n}} 
\overline{v_{j_1, j_2} v_{j_3, j_4} \cdots v_{j_{2n-1}, j_{2n}} }].
$$
Our first step in the proof of Theorem \ref{thm:main} is to prove the following.

\begin{prop} \label{prop:eq:M_unitary}
For a complex number $z$, let
\begin{equation} \label{eq:M_unitary}
M(\bm{i},\bm{j};z)=
 \sum_{\begin{subarray}{c} \sigma \in S_{2n} \\ \bm{j}=\bm{i}^\sigma \end{subarray}}
\sum_{\tau \in S_{2n}}
z^{\ell'(\tau)} \mathrm{Wg}^{U}_{2n}(\tau^{-1} \sigma;z),
\end{equation}
where $\ell'(\tau)$ is the length of the coset-type of $\tau$.
Then $M_N(\bm{i},\bm{j})=[M(\bm{i},\bm{j};z)]_{z=N}$.
\end{prop}

\begin{proof}
Equation \eqref{eq:v_u}, Lemma \ref{lem:WgFormula}, and a change of the order of sums give
\begin{align*}
M_N(\bm{i},\bm{j})=& \sum_{\bm{k} \in [N]^n} \sum_{\bm{l} \in [N]^n}
\Bigg[
\sum_{\begin{subarray}{c} \tau \in S_{2n} \\ \tilde{\bm{l}}=(\tilde{\bm{k}})^\tau \end{subarray}} 
\sum_{\begin{subarray}{c} \sigma \in S_{2n} \\ \bm{j}=\bm{i}^\sigma \end{subarray}} \mr{Wg}^{U}_{2n}(\tau^{-1} \sigma;z) \Bigg]_{z=N} \\
=& \Bigg[ \sum_{\begin{subarray}{c} \sigma \in S_{2n} \\ \bm{j}=\bm{i}^\sigma \end{subarray}} 
\sum_{\tau \in S_{2n}} \mr{Wg}^{U}_{2n}(\tau^{-1} \sigma;z)
\sum_{\bm{k} \in [N]^n} 
\sum_{\begin{subarray}{c} \bm{l} \in [N]^n \\ \tilde{\bm{l}}=(\tilde{\bm{k}})^\tau 
\end{subarray}} 1 \Bigg]_{z=N},
\end{align*}
and so the claim follows from  Lemma \ref{lem:setA} below.
\end{proof}

\begin{lem} \label{lem:setA}
Let $\tau \in S_{2n}$ and let
$A(\tau,N)=\{(\bm{k},\bm{l}) \in [N]^n \times [N]^n \ | \ 
\tilde{\bm{l}}=(\tilde{\bm{k}})^\tau  \}$.
Then $|A(\tau,N)|=N^{\ell'(\tau)}$.
\end{lem}

\begin{proof}
Given $(\bm{k},\bm{l}) \in A(\tau,N)$,
 write $\bm{k}=(k_1,k_2,\dots,k_n)$, $\bm{l}=(l_1,l_2,\dots,l_n)$, and
$\tilde{\bm{k}}=(\tilde{k}_1,\tilde{k}_2,\tilde{k}_3,\tilde{k}_4,\dots,\tilde{k}_{2n-1},\tilde{k}_{2n})$.
Then it follows immediately that
\begin{equation} \label{condition:edge}
l_j= \tilde{k}_{\tau(2j-1)}= \tilde{k}_{\tau(2j)}
\qquad \text{and} \qquad  k_j=\tilde{k}_{2j-1} 
=\tilde{k}_{2j} \qquad (j=1,2,\dots,n).
\end{equation}
In particular, $\bm{l}$ is determined by $\bm{k}$ and $\tau$.

Recall the definition of the graph $\Gamma(\tau)$ associated with $\tau$.
Given a sequence $\bm{k} \in [N]^n$,
we  assign numbers $\tilde{k}_1, \tilde{k}_2,\dots,\tilde{k}_{2n}$
to vertices $1,2,\dots,2n$ of $\Gamma(\tau)$.
Then the condition \eqref{condition:edge} implies that
two vertices $r$ and $s$
belong to the same connected component of $\Gamma(\tau)$
if and only if $\tilde{k}_r=\tilde{k}_s$.

For example, let $n=4$ and $\tau= \(\begin{smallmatrix} 1 & 2 & 3 & 4 & 5 & 6 & 7 & 8 \\
4 & 1 & 6 & 5 & 2 & 8 & 7 & 3 \end{smallmatrix}\) \in S_8$.
In this case, \eqref{condition:edge} gives the chain of {\it two} equalities
$\tilde{k}_4=\tilde{k}_1=\tilde{k}_2=\tilde{k}_8
=\tilde{k}_7=\tilde{k}_3=\tilde{k}_4$, 
$\tilde{k}_6=\tilde{k}_5=\tilde{k}_6$,
so we have $|A(\tau,N)|=|\{(k_1,k_2,k_3,k_4) \in [N]^4 \ | \ k_1=k_2=k_4\}|=N^2$. 
On the other hand, the graph $\Gamma(\tau)$ consists of {\it two} connected components,
one of which has vertices numbered by $4,1,2,8,7,3$ and the other has vertices numbered by $6,5$.

Hence, in general, $|A(\tau,N)|$ coincides with $N^{\ell'(\tau)}$, where $\ell'(\tau)$ is 
by definition the number of connected 
components of $\Gamma(\tau)$.
\end{proof}

\subsection{Expression in terms of orthogonal Weingarten functions}

Substituting \eqref{eq:ExpandN} and \eqref{eq:unitaryWg}
into \eqref{eq:M_unitary} gives
\begin{align*}
M(\bm{i},\bm{j};z)=& 
\sum_{\begin{subarray}{c} \sigma \in S_{2n} \\ \bm{j}=\bm{i}^\sigma \end{subarray}}
\sum_{\tau \in S_{2n}} \Big[
\frac{2^n n!}{(2n)!} \sum_{\lambda \vdash n}
f^{2\lambda} C'_\lambda(z) \omega^\lambda(\tau) \Big] \,
\Bigg[ \frac{1}{(2n)!} 
\sum_{\mu \vdash 2n}
\frac{f^{\mu}}{C_\mu(z)} \chi^\mu(\tau^{-1}\sigma) \Bigg]\\
=& 
\sum_{\begin{subarray}{c} \sigma \in S_{2n} \\ \bm{j}=\bm{i}^\sigma \end{subarray}}
\sum_{\lambda \vdash n} 
\sum_{\mu \vdash 2n}
\frac{2^n n!}{ \{ (2n)!\}^2}
f^{2\lambda}f^{\mu} \frac{C'_\lambda(z)}{C_\mu(z)} (\omega^\lambda * \chi^\mu)(\sigma),
\end{align*}
so Lemma \ref{lem:omega-chi} gives
$$
M(\bm{i},\bm{j};z)=
\frac{2^n n!}{(2n)!}
\sum_{\begin{subarray}{c} \sigma \in S_{2n} \\ \bm{j}=\bm{i}^\sigma \end{subarray}}
\sum_{\lambda \vdash n}
 f^{2\lambda} \frac{C_\lambda'(z)}{C_{2\lambda}(z)} \omega^\lambda(\sigma).
$$
Since $C_{2\lambda}(z)=C'_\lambda(z)C'_\lambda(z+1)$,
it follows from \eqref{eq:DefWg} that
$$
M(\bm{i},\bm{j};z)=\frac{2^n n!}{(2n)!}
\sum_{\begin{subarray}{c} \sigma \in S_{2n} \\ \bm{j}=\bm{i}^\sigma \end{subarray}}
\sum_{\lambda \vdash n }
\frac{f^{2\lambda}}{C'_\lambda(z+1)} \omega^\lambda(\sigma)
=  \sum_{\begin{subarray}{c} \sigma \in S_{2n} \\ \bm{j}=\bm{i}^\sigma \end{subarray}}
\mathrm{Wg}^{O}_{n}(\sigma;z+1).
$$
Hence, the last equation together with Proposition \ref{prop:eq:M_unitary} implies Theorem \ref{thm:main}.

\section{Moments of a single matrix element} \label{section:single moment}

We next obtain explicit values for the moment
of a single matrix element $\bE[|v_{ij}|^{2n}]$.

\begin{thm} \label{thm:diagonal}
Let $N \ge 1$ and let $v_{ii}$ be a diagonal entry of an $N \times N$ COE matrix. Then,
for $n \ge 1$,
$$
\bE[|v_{ii}|^{2n}] = \frac{2^n n!}{(N+1)(N+3) \cdots (N+2n-1)}.
$$
\end{thm}

\begin{proof}
Theorem \ref{thm:main} implies 
$$
\bE[|v_{ii}|^{2n}]
=M_N((\underbrace{i,\dots,i}_{2n}),(\underbrace{i,\dots,i}_{2n})) = \Bigg[\sum_{\sigma \in S_{2n}}
 \mathrm{Wg}^{O}_{n}(\sigma;z) \Bigg]_{z=N+1}.
$$
Since $\chi^{(2n)} \equiv 1$,  
Lemma \ref{lem:omega-chi} implies that
$\sum_{\sigma \in S_{2n}} \omega^\lambda(\sigma)
=(\omega^\lambda *\chi^{(2n)}) (\mathrm{id}_{2n})= 
(2n)! \delta_{\lambda, (n)}$,
and so we have
\begin{align*}
&\sum_{\sigma \in S_{2n}} \mathrm{Wg}^{O}_{n}(\sigma;z) =
\frac{2^n n!}{(2n)!}
\sum_{\lambda \vdash n}
\frac{f^{2\lambda}}{C'_\lambda(z)} \sum_{\sigma \in S_{2n}}\omega^\lambda(\sigma) \\
=& \frac{ 2 ^n n! f^{(2n)}}{C'_{(n)}(z)} 
=\frac{2^n n!}{z(z+2) \cdots (z+2(n-1))},
\end{align*}
which implies the desired result.
\end{proof}

\begin{thm} \label{thm:off-diagonal}
Let $N\ge 2$ and let $v_{ij}$ be an off-diagonal entry of an $N \times N$ COE matrix. Then,
for $n \ge 1$,
$$
\bE[|v_{ij}|^{2n}] =\frac{n!}{(N+2n-1) \prod_{k=0}^{n-2}(N+k)}.
$$
\end{thm}

The proof of Theorem \ref{thm:off-diagonal} is more difficult than 
that of Theorem \ref{thm:diagonal}.
We first show the following lemma.

\begin{lem} 
Under the same assumptions as in Theorem \ref{thm:off-diagonal},
\begin{equation} \label{eq:off-diagonal_lemma}
\bE[|v_{ij}|^{2n}] = 
\frac{(n!)^2}{(2n)!} \sum_{
\begin{subarray}{c} \lambda \vdash n \\ \ell(\lambda) \le 2
\end{subarray}} f^{2\lambda} \frac{C_\lambda'(2)}{C'_\lambda(N+1)}.
\end{equation}
\end{lem}

\begin{proof}
Let $\bm{i}=(i,j,i,j,\dots,i,j) \in [N]^{2n}$, where $i \not=j$.
If $\sigma \in S_{2n}$ satisfies $\bm{i}^\sigma =\bm{i}$, then
it permutates odd numbers and even numbers, and so 
it is uniquely expressed as
$$
\sigma(2k-1)= 2 \pi_1(k)-1, \qquad \sigma(2k)= 2 \pi_2(k), \qquad k=1,2,\dots,n
$$
for some $\pi_1,\pi_2 \in S_{n}$.
In this case,
it is not difficult to see that
the coset-type of $\sigma$ coincides with the cycle-type of $\pi_1 \pi_2^{-1}$.
It follows from Theorem \ref{thm:main} that 
$\bE[|v_{ij}|^{2n}] = [\mcal{M}(\bm{i},\bm{i};z)]_{z=N+1}$
with
\begin{align*}
\mcal{M}(\bm{i},\bm{i};z)=& \sum_{\begin{subarray}{c} \sigma \in S_{2n} \\
\bm{i}^\sigma=\bm{i} \end{subarray}} 
\mathrm{Wg}^O_n(\sigma;z) 
= \sum_{\pi_1, \pi_2 \in S_n} \mathrm{Wg}^{O}_n (\mathrm{Cycle}(\pi_1 \pi_2^{-1});z) \\
=& n! \sum_{\pi \in S_n} \mathrm{Wg}^{O}_n (\mathrm{Cycle}(\pi);z),
\end{align*}
where $\mathrm{Cycle}(\pi)$ is the cycle-type of $\pi$
and $\mathrm{Wg}^{O}_n(\mu;z)$ is the value of $\mathrm{Wg}^O_n(\sigma;z)$ at $\sigma \in H_\mu$.
The number of permutations in $S_n$ of cycle-type $\mu$ is $n!/z_\mu$; hence,
$$
\mcal{M}(\bm{i},\bm{i};z) = n! \sum_{\mu \vdash n} \frac{n!}{z_\mu} \mathrm{Wg}^{O}_n(\mu;z).
$$
Since $|H_\mu|=(2^n n!)^2 /(2^{\ell(\mu)} z_\mu)$, we have
$$
\mcal{M}(\bm{i},\bm{i};z)  = 2^{-4n} \sum_{\mu \vdash n} 2^{\ell(\mu)} |H_\mu| \mathrm{Wg}^{O}_n
(\mu;z) = 2^{-4n} \sum_{\sigma \in S_{2n}} 2^{\ell'(\sigma)} \mathrm{Wg}^{O}_n(\sigma;z)
$$
by the decomposition $S_{2n}= \bigsqcup_{\mu \vdash n} H_\mu$.
Using \eqref{eq:ExpandN}, \eqref{eq:DefWg}, and Lemma \ref{lem:omega-chi} gives
\begin{align*}
\mcal{M}(\bm{i},\bm{i};z)=& 
\(\frac{n!}{(2n)!} \)^2 
\sum_{\lambda \vdash n} 
\sum_{\mu \vdash n }
f^{2\lambda} f^{2\mu} \frac{C_\lambda'(2)}{C'_\mu(z)}
(\omega^\lambda *\omega^\mu)(\mathrm{id}_{2n}) \\
=&\frac{(n!)^2}{(2n)!} \sum_{\lambda \vdash n } 
f^{2\lambda} \frac{C_\lambda'(2)}{C'_\lambda(z)},
\end{align*}
where we used $\omega^\lambda(\mathrm{id}_{2n})=1$ to obtain the last equality.
Since $C_\lambda'(2) =0$ if $\ell(\lambda) >2$,
the range of the last sum can be restricted to $\{\lambda \vdash n \ | \ \ell(\lambda) \le 2\}$.
Since $C'_\lambda(N+1)$ is nonzero for every $\lambda$ in the restricted set,
we have
$$
\bE[|v_{ij}|^{2n}] = [\mcal{M}(\bm{i},\bm{i};z)]_{z=N+1} = 
\frac{(n!)^2}{(2n)!} \sum_{
\begin{subarray}{c} \lambda \vdash n \\ \ell(\lambda) \le 2
\end{subarray}} f^{2\lambda} \frac{C_\lambda'(2)}{C'_\lambda(N+1)}.
$$
\end{proof}

We next evaluate the sum on the right-hand side of \eqref{eq:off-diagonal_lemma}.

If $\mu$ is a partition of $n$,
the well-known hook-length formula states that $f^{\mu} = n! /h(\mu)$,
where $h(\mu)$ is the hook-length product
$$
h(\mu)= \prod_{(i,j) \in \mu} (\mu_i+\mu_j'-i-j+1)
$$
and  $\mu'=(\mu_1',\mu_2',\dots)$ is the conjugate partition of $\mu$.

A partition $\lambda$ of $n$ with $\ell(\lambda)\le 2$
can be expressed as $\lambda=(n-r,r)$ with some $0 \le r \le [n/2]$,
where $[n/2]$ is the integer part of $n/2$.
We now fix $n, N$ and 
let
$$
A_r= \frac{f^{(2n-2r,2r)} C_{(n-r,r)}'(2)}{(2n)! \cdot C_{(n-r,r)}'(N+1)}
=\frac{C_{(n-r,r)}'(2)}{h((2n-2r,2r)) C_{(n-r,r)}'(N+1)}
$$
for each $0 \le r \le [n/2]$.
A direct computation gives
$$
A_r= \frac{2n-4r+1}{(2r)!! \cdot (2n-2r+1)!!} 
\prod_{j=1}^{n-r} \frac{1}{N+2j-1} \cdot \prod_{j=0}^{r-1} \frac{1}{N+2j}.
$$
Here, 
$(2k-1)!!= \prod_{j=1}^k (2j-1)$ and $(2k)!! = \prod_{j=1}^k (2j)$ with $0!!=1$.

\begin{lem} \label{lem:sumA}
For each $0 \le k \le [n/2]$,
\begin{align*}
A_0+A_1+\cdots+A_k=& \frac{1}{(2k)!! \cdot (2(n-k)-1)!!}\\
& \times \frac{1}{(N+2n-1)} \cdot \prod_{j=1}^{n-k-1} \frac{1}{N+2j-1}  \cdot 
\prod_{j=0}^{k-1} \frac{1}{N+2j}.
\end{align*}
\end{lem}

\begin{proof}
We show the lemma by induction on $k$.
Observe that the case at $k=0$ holds true:
$$
A_0 = \frac{1}{(2n-1)!!}\prod_{j=1}^n  \frac{1}{N+2j-1}. 
$$
We next consider the identity at $k+1$.
From the induction assumption, we have
\begin{align*}
&A_0+A_1+\cdots+A_k+A_{k+1} \\
=& 
\frac{1}{(2k)!! \cdot (2n-2k-1)!!} \cdot 
\frac{1}{(N+2n-1)} \cdot \prod_{j=1}^{n-k-1} \frac{1}{N+2j-1}  \cdot 
\prod_{j=0}^{k-1} \frac{1}{N+2j} \\
&+ \frac{2n-4k-3}{(2k+2)!! \cdot (2n-2k-1)!!} 
\prod_{j=1}^{n-k-1} \frac{1}{N+2j-1} \cdot \prod_{j=0}^{k} \frac{1}{N+2j}.
\end{align*}
By a reduction, 
the last equation is 
\begin{align*}
=& \frac{1}{(2k+2)!! (2n-2k-1)!!} \frac{ (2k+2)(N+2k) + (2n-4k-3)(N+2n-1)  
}
{ (N+2n-1)\prod_{j=1}^{n-k-1} (N+2j-1)  \cdot \prod_{j=0}^{k} (N+2j)}\\
=&\frac{1}{(2k+2)!! (2n-2k-1)!!} \frac{(2n-2k-1)(N+2n-2k-3)}{ (N+2n-1)\prod_{j=1}^{n-k-1} (N+2j-1)  \cdot \prod_{j=0}^{k} (N+2j)} \\
=& \frac{1}{(2k+2)!! (2n-2k-3)!!} \frac{1}{ (N+2n-1)\prod_{j=1}^{n-k-2} (N+2j-1)  \cdot \prod_{j=0}^{k} (N+2j)},
\end{align*}
which concludes the desired identity.
\end{proof}

Finally, we will prove Theorem \ref{thm:off-diagonal}.

\begin{proof}[Proof of Theorem \ref{thm:off-diagonal}]
From \eqref{eq:off-diagonal_lemma} and the definition of $A_r$, we have
$$
\bE[|v_{ij}|^{2n}]= (n!)^2 \sum_{r=0}^{[n/2]} A_r.
$$
We first suppose that $n$ is even: $n=2m$. Using Lemma \ref{lem:sumA}, we have
\begin{align*}
\sum_{r=0}^{m} A_r =& \frac{1}{(2m)!! \cdot (2m-1)!!}
\frac{1}{(N+2n-1)} \cdot \prod_{j=1}^{m-1} \frac{1}{N+2j-1}  \cdot 
\prod_{j=0}^{m-1} \frac{1}{N+2j} \\
=& \frac{1}{(2m)!} \frac{1}{(N+2n-1)} \cdot \prod_{k=0}^{2m-2} \frac{1}{N+k},
\end{align*}
which implies Theorem \ref{thm:off-diagonal}.
The same result holds for the odd $n$ case by a similar calculation.
Thus, the proof of Theorem \ref{thm:off-diagonal} is complete.
\end{proof}

\begin{remark}
Let $o_{ii}$ and $o_{jj}$ $(i \not=j)$ be two diagonal entries of 
a Haar-distributed orthogonal matrix from $O(N)$.
Theorem 5.1 in \cite{CM} gives that
$$
\bE[ (o_{ii}+o_{jj})^{2n}] = 
\sum_{
\begin{subarray}{c} \lambda \vdash n \\ \ell(\lambda) \le 2
\end{subarray}} f^{2\lambda} \frac{C_\lambda'(2)}{C'_\lambda(N)}
$$
for all $n \ge 1$ and $N \ge 2$.
From the discussion in this section, we have obtained a new identity
$$
\bE[ (o_{ii}+o_{jj})^{2n}]= \frac{(2n)!}{n!} \cdot
\frac{1}{(N+2n-2) \prod_{k=0}^{n-2}(N+k-1)}.
$$
\end{remark}

\section{Asymptotic behavior} \label{section:asymptotic}

We now give an asymptotic expansion for $M_N(\bm{i},\bm{j})$ for the case $N \to \infty$.
In this section,
do not confuse Landou's big $O$ notation $O(N^{-k})$ with the symbol for
the orthogonal group $O(N)$.

As we saw in subsection \ref{subsec:Wg},
$$
\mathrm{Wg}^{O(N)}_n(\sigma)=\mathrm{Wg}^{O}_n(\sigma;N)=
\frac{2^n n!}{(2n)!} 
\sum_{\lambda \vdash n}
\frac{f^{2\lambda}}{C_\lambda'(N)} 
\omega^\lambda(\sigma) \qquad (\sigma \in S_{2n})
$$
if $N$ is sufficiently large.
We already know the asymptotic behavior of $\mathrm{Wg}^{O(N)}_n$ as $N \to \infty$
as follows.

\begin{lem}[\cite{Mat_JMortho}]
Fix $n$ and $\sigma \in S_{2n}$. As $N \to \infty$,
\begin{enumerate}
\item If $\sigma \in H_n$, then $\mathrm{Wg}^{O(N)}_n(\sigma)= 
N^{-n} +O(N^{-n-2})$.
\item If $\sigma \in H_{(2,1^{n-2})}$, then
$\mathrm{Wg}^{O(N)}_n (\sigma) = -N^{-n-1} +O(N^{-n-2})$.
\item Otherwise, $\mathrm{Wg}^{O(N)}_n (\sigma) = O(N^{-n-2})$.
\end{enumerate}
\end{lem}

\begin{thm} \label{thm:asym}
Fix a positive integer $n$ and fix
two sequences $\bm{i}=(i_1,\dots,i_{2n})$ and $\bm{j}=(j_1,\dots,j_{2n})$.
As $N \to \infty$,
$$
M_N(\bm{i},\bm{j}) = s(\bm{i},\bm{j}) N^{-n}-(s'(\bm{i},\bm{j})+ s(\bm{i},\bm{j})n) N^{-n-1}
+O(N^{-n-2}),
$$
where 
$$
s(\bm{i},\bm{j})= |\{\sigma \in H_{n} \ | \ \bm{j}=\bm{i}^\sigma\}|
\qquad \text{and} \qquad 
s'(\bm{i},\bm{j})= |\{\sigma \in H_{(2,1^{n-2})} \ | \ \bm{j}=\bm{i}^\sigma\}|.
$$
\end{thm}

\begin{proof}
Since $(N+1)^{-n}= N^{-n} - n N^{-n-1} +O(N^{-n-2})$,
the last lemma gives
$$
\mathrm{Wg}^{O(N+1)}_{n}(\sigma)=
\begin{cases}
N^{-n}-n N^{-n-1} + O(N^{-n-2}) & \text{if $\sigma \in H_n$}, \\
-N^{-n-1} +O(N^{-n-2}) & \text{if $\sigma \in H_{(2,1^{n-2})}$}, \\
O(N^{-n-2}) & \text{otherwise}.
\end{cases}
$$
It follows from Theorem \ref{thm:main} and the double coset decomposition
$S_{2n}=\bigsqcup_{\mu \vdash n} H_\mu$
that $M_N(\bm{i}, \bm{j})$ is 
\begin{align*}
& \sum_{\begin{subarray}{c} \sigma \in H_{n} \\
\bm{j}=\bm{i}^\sigma \end{subarray}}
\mathrm{Wg}^{O(N+1)}_n (\sigma) 
+\sum_{\begin{subarray}{c} \sigma \in H_{(2,1^{n-2})} \\
\bm{j}=\bm{i}^\sigma \end{subarray}}
\mathrm{Wg}^{O(N+1)}_n (\sigma)
+ \sum_{\begin{subarray}{c} \mu \vdash n \\
\mu \not= (1^n), \ (2,1^{n-2}) \end{subarray}}
\sum_{\begin{subarray}{c} \sigma \in H_\mu \\
\bm{j}=\bm{i}^\sigma \end{subarray}}
\mathrm{Wg}^{O(N+1)}_n (\sigma)
\\
=& 
s(\bm{i},\bm{j}) (N^{-n} -n N^{-n-1}) - s'(\bm{i},\bm{j}) N^{-n-1} +O(N^{-n-2})
\end{align*}
as $N \to \infty$.
\end{proof}

If $s(\bm{i},\bm{j})>0$, then there exists a $\zeta \in H_n$ such that $\bm{j}=\bm{i}^\zeta$.
Fixing $\zeta$, the correspondence $\iota:\sigma \mapsto \sigma \zeta$ defines 
a bijection 
$\iota: \{\sigma \in H_n  \ | \ \bm{i}=\bm{i}^\sigma \} \to \{\tau \in H_n \ | \ 
\bm{j}=\bm{i}^\tau\}$.
Therefore, in this case, $s(\bm{i},\bm{j})$ coincides with
the order of the stabilizer subgroup of $\bm{i}$ in $H_n$.

\begin{example} \label{ex:asymM}
\begin{enumerate}
\item
Let $\bm{i}=(i,i,\dots,i) \in [N]^{2n}$.
Then $s(\bm{i},\bm{i})=|H_n|=2^n n!$ and
$s'(\bm{i},\bm{i})=|H_{(2,1^{n-2})}| = n(n-1)2^n n!$,
and so we have
$$
\bE[|v_{ii}|^2]= 2^n n! (N^{-n}- n^2 N^{-n-1}) +O(N^{-n-2}).
$$
\item Let $ i \not=j$ and let
$\bm{i}= (i,j,i,j,\dots,i,j)$.
Then
$s(\bm{i},\bm{i})= n!$ and $s'(\bm{i},\bm{i})= \binom{n}{2}  n!$.
Therefore,
$$
\bE[|v_{ij}|^{2n}]= n! \( N^{-n} -\frac{n(n+1)}{2} N^{-n-1} \) +O(N^{-n-2}).
$$
\item Let $\bm{i}=(i_1,i_2,\dots,i_{2n})$ be a sequence
with distinct numbers. Then $s(\bm{i},\bm{i})=1$ and $s'(\bm{i},\bm{i})=0$; hence,
$$
\bE [| v_{i_1 i_2} v_{i_3 i_4} \cdots v_{i_{2n-1},i_{2n}}|^{2}]
= N^{-n} -n N^{-n-1} +O(N^{-n-2}).
$$
\item Let $\bm{i}=(i,j,i,j)$ and $\bm{j}=(i,i,j,j)$ with $i \not=j$. 
Then $s(\bm{i},\bm{j})=0$ and $s'(\bm{i},\bm{j})=4$, and so we have
$\bE[v_{ij}^2 \overline{v_{ii} v_{jj}}] = -4 N^{-3} +O(N^{-4})$.
\end{enumerate}
\end{example}

\medskip

A standard complex Gaussian random variable $Z$ is defined by
$Z=(X_1 + \sqrt{-1} X_2)/\sqrt{2}$,
where $X_1$ and $X_2$ are independent and standard (real) Gaussian random variables.
Then $Z$ satisfies $\bE[ Z^n \overline{Z}^m]= \delta_{nm} n!$ for any positive integers $n$ and $m$.
Lemma \ref{lem:COEcondition} implies 
$\bE[v_{ii}^n \overline{v_{ii}}^m]=\bE[v_{ij}^n \overline{v_{ij}}^m]=0$
with $i \not=j$
if $n \not=m$,
and the first two claims in Example \ref{ex:asymM} 
(or Theorem \ref{thm:diagonal}, \ref{thm:off-diagonal})
imply
$$
\lim_{N \to \infty} \(\frac{N}{2}\)^n \bE[ |v_{ii}|^{2n}] =
\lim_{N \to \infty} N^n \bE[ |v_{ij}|^{2n}] 
= n!.
$$ 
Hence, by the Carleman criterion for moment problems (see the proof of 
Theorem 3.1 in \cite{PR}),
we obtain an algebraic proof of the following result,
which was proved by Jiang \cite{Jiang}.

\begin{cor}
Let $V^{(N)}=(v_{ij}^{(N)})$, $N \ge 1$, 
be a sequence of COE matrices. 
Fix positive integers $i \not= j$. 
As $N \to \infty$, 
the random variables 
$\sqrt{N/2} v_{ii}^{(N)}$ and $\sqrt{N} v_{ij}^{(N)}$ 
converge weakly to a standard complex Gaussian random variable.
\end{cor}

\section{Example} \label{section:example}

Finally, we give some examples for small $n$ by applying Theorem \ref{thm:main}.
Let $V=(v_{ij})$ be an $N \times N$ COE matrix as usual.
We have already obtained explicit values of moments of the form
$\bE[|v_{ij}|^{2n}]$ when describing Theorems \ref{thm:diagonal} and \ref{thm:off-diagonal},
and so let us consider other examples.

Let $\bm{i}=(i_1,i_2,\dots,i_{2n})$ and $\bm{j}=(j_1,j_2,\dots,j_{2n})$ be
sequences in $[N]^{2n}$.
Let 
$$
s_\mu(\bm{i},\bm{j}) =|\{\sigma \in H_\mu \ | \ \bm{j}=\bm{i}^\sigma\}|
$$
for each $\mu \vdash n$.
Note that $s_{(1^n)}(\bm{i},\bm{j}) = s(\bm{i},\bm{j})$ and 
$s_{(2, 1^{n-2})}(\bm{i},\bm{j}) =s'(\bm{i},\bm{j})$ in the notation of 
Theorem \ref{thm:asym}.
Then Theorem \ref{thm:main} implies 
$$
M_N(\bm{i},\bm{j})=\sum_{\mu \vdash n}
s_\mu(\bm{i},\bm{j}) \mathrm{Wg}^{O}_n(\mu;N+1),
$$
which is an identity in $N$, and $N$ is regarded as 
an indeterminate rather than a fixed positive integer here. 
Some explicit values for $\mathrm{Wg}^{O}_n(\mu;z)$ can be seen in 
Example \ref{ex:WgO} or in \cite{CS,CM}.

\begin{example} (degree 1.)
If $i_1,i_2,j_1,j_2$ are positive integers in $[N]$,
$$
\bE[v_{i_1 i_2} \overline{v_{j_1 j_2}}]= M_N((i_1,i_2),(j_1,j_2))=
\begin{cases}
\frac{2}{N+1} & \text{if $i_1=i_2=j_1=j_2$},\\
\frac{1}{N+1} & \text{if $i_1 \not=i_2$, $j_1 \not=j_2$, and $\{i_1,i_2\}=\{j_1,j_2\}$},\\
0 & \text{otherwise}.
\end{cases}
$$
\end{example}

\begin{example} 
(degree 2.) 
Let $i \not=j$.
Put
\begin{align*}
L_1=& \left\{ \(\begin{smallmatrix} 1 & 2 & 3 & 4 \\
1 & 2 & 3 & 4 \end{smallmatrix}\),\(\begin{smallmatrix} 1 & 2 & 3 & 4 \\
2 & 1 & 3 & 4 \end{smallmatrix}\),
\(\begin{smallmatrix} 1 & 2 & 3 & 4 \\
1 & 2 & 4 & 3 \end{smallmatrix}\),
\(\begin{smallmatrix} 1 & 2 & 3 & 4 \\
2 & 1 & 4 & 3 \end{smallmatrix}\)
\right\}, \\
L_2=& \left\{ \(\begin{smallmatrix} 1 & 2 & 3 & 4 \\
1 & 3 & 2 & 4 \end{smallmatrix}\),\(\begin{smallmatrix} 1 & 2 & 3 & 4 \\
3 & 1 & 2 & 4 \end{smallmatrix}\),
\(\begin{smallmatrix} 1 & 2 & 3 & 4 \\
1 & 3 & 4 & 2 \end{smallmatrix}\),
\(\begin{smallmatrix} 1 & 2 & 3 & 4 \\
3 & 1 & 4 & 2 \end{smallmatrix}\)
\right\}.
\end{align*}
Using the list of $\mathrm{Wg}^O_2(\mu;z)$ given as Example \ref{ex:WgO},
we can compute
\begin{align*}
\bE[|v_{ii} v_{jj}|^2] = &M_N((i,i,j,j),(i,i,j,j))
= \sum_{\sigma \in L_1} \mathrm{Wg}_2^{O}(\sigma;N+1)  \\
=& 4 \cdot \mathrm{Wg}^O_2([1^2];N+1) 
= \frac{4(N+2)}{(N+1)N(N+3)},
\end{align*}
\begin{align*}
\bE[v_{ij}^2 \overline{v_{ii} v_{jj}}] =& M_N((1,2,1,2),(1,1,2,2))
= \sum_{\sigma \in L_2}\mathrm{Wg}_2^{O}(\sigma;N+1) 
\\
=&4 \cdot \mathrm{Wg}^O_2([2];N+1) 
= \frac{-4}{(N+1)N(N+3)}.
\end{align*}
\end{example}

\begin{example}
(degree 3.)
Using Theorem \ref{thm:main} and the list of $\mathrm{Wg}^O_3(\mu;z)$ given as Example \ref{ex:WgO},
we have
$$
\bE [|v_{11} v_{22} v_{33}|^2]= 8 \cdot \mathrm{Wg}^{O}_3 ([1^3];N+1) =
\frac{8(N^2+5N+2)}{(N-1)N(N+1)(N+3)(N+5)},
$$
\begin{align*}
\bE[ v_{12}v_{34}v_{56} \overline{v_{13}v_{26}v_{45}}] 
=& \mathrm{Wg}^{O}_3 \( \( \begin{smallmatrix}
1 & 2 & 3 & 4 & 5 & 6 \\
1 & 3 & 2 & 6 & 4 & 5 
\end{smallmatrix}\);N+1\)
= \mathrm{Wg}^{O}_3 ([3];N+1) \\
=&\frac{2}{(N-1)N(N+1)(N+3)(N+5)},
\end{align*}
\begin{align*}
 \bE[ |v_{11}^2 v_{22}|^2]  
=& 16 \cdot \mathrm{Wg}^{O}_3 ([1^3];N+1)
+ 32 \cdot \mathrm{Wg}^{O}_3 ([2,1];N+1) \\
=& \frac{16(N+4)}{N(N+1)(N+3)(N+5)}. 
\end{align*}
\end{example}


\end{document}